\documentclass[a4paper,leqno]{article}
\usepackage[dvips]{graphicx}
\usepackage{amssymb}
\usepackage{amsmath}
\usepackage{mathrsfs}
\usepackage{color}
\usepackage{theorem}
\usepackage{comment}
\usepackage{url}
\usepackage{ascmac}
\newcommand{\pref}[1]{(\ref{#1})}
\newcommand{\be}{\begin{equation}}
\newcommand{\ee}{\end{equation}}

\newcommand{\vecf}{\mbox{\boldmath $ f $}}

\newcommand{\vecu}{\mbox{\boldmath $ u $}}

\newcommand{\vectau}{\mbox{\boldmath $ \tau $}}

\renewcommand{\epsilon}{\varepsilon}

\newtheorem{thm}{Theorem}[section]

\newtheorem{cor}{Corollary}[section]

\theorembodyfont{\rmfamily}

\title{Cosine formula for generalized O'Hara's energies}
\author{Takeyuki Nagasawa
\thanks{Saitama University,
Japan.
The author is supported by Grant-in-Aid for Scientific Research (C) (No.17K05310),
Japan Society for the Promotion of Science.
}
}
\date{}
\begin{document}
\maketitle
\begin{abstract}
In this short article,
we extend the cosine formula for the M\"{o}bius energy to generalized O'Hara energies.
The newly derived formula gives us a condition for which the right circle minimizes the energy under the length-constraint.
Furthermore,
it shows us how far the energy is from the M\"{o}bius invariant property.
\\
keywords:
O'Hara's energy,
M\"{o}bius energy,
knot energy,
cosine formula 
\end{abstract}
\section{Introduction}
\par
In a series of papers \cite{OH1}--\cite{OH3},
O'Hara proposed the energy functional
\[
	\mathcal{E}_{( \alpha, p ) } ( \vecf )
	=
	\iint_{ ( \mathbb{R} / \mathcal{L} \mathbb{Z} )^2 }
	\left(
	\frac 1 { \| \vecf ( s_1 ) - \vecf ( s_2 ) \|_{ \mathbb{R}^3 }^\alpha }
	-
	\frac 1 { \mathscr{D} ( \vecf ( s_1 ) , \vecf ( s_2 ) )^\alpha }
	\right)^p
	d s_1 d s_2
\]
for an $ \mathbb{R}^3 $-valued function $ \vecf $ on $ \mathbb{R} / \mathcal{L} \mathbb{Z} $ representing a knot parametrized by the arc-length.
Among these,
the energy $ \mathcal{E}_{ ( 2,1 ) } $ is prominent for the invariance under M\"{o}bius transformations.
This property had been first identified by Freedman-He-Wang \cite{FHW},
while Doyle and Schramm later gave an alternative proof showing the expression
\be
	\mathcal{E}_{( 2,1 ) } ( \vecf )
	=
	\iint_{ ( \mathbb{R} / \mathcal{L} \mathbb{Z} )^2 }
	\frac { 1 - \cos \varphi ( s_1 , s_2 ) } { \| \vecf ( s_1 ) - \vecf ( s_2 ) \|_{ \mathbb{R}^n }^2 } \, d s_1 d s_2 + 4
	,
	\label{cosine_formula_of_Moebius_energy}
\ee
where $ \varphi $ is the conformal angle.
This is called the {\it cosine formula},
and was presented in \cite{KS}.
The conformal angle is defined as follows.
Let $ C_{12} $ be the circle contacting a knot $ \mathrm{Im} \vecf $ at $ \vecf ( s_1 ) $ and passing through $ \vecf ( s_2 ) $.
We define the circle $ C_{21} $ similarly.
The angle $ \varphi ( s_1 , s_2 ) $ is that between these two circles at $ \vecf ( s_1 ) $ (and also at $ \vecf ( s_2 ) $).
Note that it is determined from $ \vecf ( s_1 ) - \vecf ( s_2 ) $,
$ \vecf^\prime ( s_1 ) $,
and $ \vecf^\prime ( s_2 ) $.
It is clearly invariant under M\"{o}bius transformations.
The invariance of $ \displaystyle{ \frac { d s_1 d s_2 } { \| \vecf ( s_1 ) - \vecf ( s_2 ) \|_{ \mathbb{R}^3 }^2 } } $ under M\"{o}bius transformations follows from that of the cross ratio.
Hence,
we can easily read the M\"{o}bius invariance property of $ \mathcal{E}_{ (2,1) } $ from \pref{cosine_formula_of_Moebius_energy}.
\par
In this article,
we discuss an analogue for generalized O'Hara energies.
The motivation is as follows.
In the following, 
we use the notation $ \Delta $ to mean $ \Delta \vecu = \vecu ( s_1 ) - \vecu ( s_2 ) $,
and $ \| \cdot \| = \| \cdot \|_{ \mathbb{R}^n } $,
$ \langle \cdot , \cdot \rangle = \langle \cdot , \cdot \rangle_{ { \bigwedge \! }^2 \, \mathbb{R}^n } $.
According to \cite{Blatt2},
if $ \mathcal{E}_{(2,1)} ( \vecf ) < \infty $,
then the unit tangent vector $ \vectau = \vecf^\prime $ exists almost everywhere.
The author,
together with Ishizeki \cite{IshizekiNagasawaI,IshizekiNagasawaIII} showed a decomposition of $ \mathcal{E}_{(2,1)} $ for an $ \mathbb{R}^n $-valued function:
\begin{align}
	\mathcal{E}_{(2,1)} ( \vecf ) 
	= & \
	\mathcal{E}_{(2,1),1} ( \vecf ) + \mathcal{E}_{(2,1),2} ( \vecf )
	+
	4
	,
	\label{docomposition_of_Moebius_energy}
	\\
	\mathcal{E}_{(2,1),1} ( \vecf )
	= & \
	\iint_{ ( \mathbb{R} / \mathcal{L} \mathbb{Z} )^2 }
	\frac
	{ \| \Delta \vectau \|^2 }
	{ 2 \| \Delta \vecf \|^2 }
	\,
	d s_1 d s_2 ,
	\nonumber
	\\
	\mathcal{E}_{(2,1),2} ( \vecf )
	= & \
	\iint_{ ( \mathbb{R} / \mathcal{L} \mathbb{Z} )^2 }
	\frac 2
	{ \| \Delta \vecf \|_{ \mathbb{R}^n }^2 }
	\left\langle
	\vectau ( s_1 ) \wedge \frac { \Delta \vecf } { \| \Delta \vecf \| }
	,
	\vectau ( s_2 ) \wedge \frac { \Delta \vecf } { \| \Delta \vecf \| }
	\right\rangle
	d s_1 d s_2
	.
	\nonumber
\end{align}
Each energy $ \mathcal{E}_{(2,1),i} $ is also M\"{o}bius invariant.
The decomposition was extended for the generalized O'Hara energy
\[
	\mathcal{E}_\Phi ( \vecf )
	=
	\iint_{ ( \mathbb{R} / \mathcal{L} \mathbb{Z} )^2 }
	\left(
	\frac 1 { \Phi ( \| \Delta \vecf \| ) }
	-
	\frac 1 { \Phi ( \mathscr{D} ( \vecf ( s_1 ) , \vecf ( s_2 ) ) ) }
	\right)
	d s_1 d s_2
	,
\]
under suitable assumptions on the function $ \Phi $ from $ \mathbb{R}_+ = \{ x \in \mathbb{R} \, | \, x > 0 \} $ to itself,
as
\begin{align}
	\mathcal{E}_\Phi ( \vecf ) 
	= & \
	\mathcal{E}_{\Phi,1 } ( \vecf ) + \mathcal{E}_{\Phi,2} ( \vecf )
	+
	2 \mathcal{L} \int_{ \frac { \mathcal{L} } 2 }^\infty \frac { dt } { \Phi (t) }
	,
	\label{decomposition_of_generalized_O'Hara_energy}
	\\
	\mathcal{E}_{\Phi,1} ( \vecf )
	= & \
	\iint_{ ( \mathbb{R} / \mathcal{L} \mathbb{Z} )^2 }
	\frac
	{ \| \Delta \vectau \|^2 }
	{ 2 \Phi ( \| \Delta \vecf \| ) }
	\,
	d s_1 d s_2 ,
	\nonumber
	\\
	\mathcal{E}_{\Phi,2} ( \vecf )
	= & \
	\iint_{ ( \mathbb{R} / \mathcal{L} \mathbb{Z} )^2 }
	\left(
	\frac 1 { \Phi ( \| \Delta \vecf \| ) }
	-
	\Lambda ( \| \Delta \vecf \| )
	\right)
	\nonumber
	\\
	& \qquad
	\times
	\left\langle
	\vectau ( s_1 ) \wedge \frac { \Delta \vecf } { \| \Delta \vecf \| }
	,
	\vectau ( s_2 ) \wedge \frac { \Delta \vecf } { \| \Delta \vecf \| }
	\right\rangle
	d s_1 d s_2
	,
	\nonumber
	\\
	\Lambda ( t )
	= & \
	- \frac 1t \int_t^\infty \frac { dx } { \Phi (x) }
	\nonumber
\end{align}
see \cite{IshizekiNagasawaV}.
Since $ \mathcal{E}_\Phi $ does not generally have the M\"{o}bius invariant property,
$ \mathcal{E}_{ \Phi , i } $ also does not.
We,
however,
can derive the first and second variational formulae from \pref{decomposition_of_generalized_O'Hara_energy} systematically as in \cite{IshizekiNagasawaII}.
The part ``$ + \, 4 $'' is common for both \pref{cosine_formula_of_Moebius_energy} and \pref{docomposition_of_Moebius_energy}.
Then \pref{decomposition_of_generalized_O'Hara_energy} suggests that the energy $ \mathcal{E}_\Phi $ also has the expression
\[
	\mathcal{E}_\Phi ( \vecf )
	=
	\iint_{ ( \mathbb{R} / \mathcal{L} \mathbb{Z} )^2 }
	\frac { 1 - \cos \psi ( s_1 , s_2 ) } { \Phi ( \| \Delta \vecf \| ) } \, d s_1 d s_2
	+
	2 \mathcal{L} \int_{ \frac { \mathcal{L} } 2 }^\infty \frac { dt } { \Phi (t) }
\]
for some angle $ \psi $ determined from $ \Delta \vecf $,
$ \vectau ( s_1 ) $,
and $ \vectau ( s_2 ) $.
The answer is affirmative,
and the formula indicates how far the energy is from the M\"{o}bius invariant.
In next section,
we will give the formula and discuss its consequences.
The proof of the formula will be given in the last section.
\section{Cosine formula}
\par
In our previous paper \cite{IshizekiNagasawaV} we proved \pref{decomposition_of_generalized_O'Hara_energy} under the assumption that
\begin{itemize}
\item[{\rm (A.1)}] $ \Phi $ is monotonically increasing. 
\item[{\rm (A.2)}]
	$ \displaystyle{ \int_x^\infty \frac { dt } { \Phi (t) } < \infty } $ for $ x > 0 $.
\item[{\rm (A.3)}]
	The function space $ W_\Phi $ is defined by
	\[
		W_\Phi
		=
		\left\{ \vecu \in W^{1,2} ( \mathbb{R} / \mathcal{L} \mathbb{Z} )
		\, \left| \,
		\iint_{ ( \mathbb{R} / \mathcal{L} \mathbb{Z} )^2 }
		\frac { \| \Delta \vecu^\prime \|_{ \mathbb{R}^n }^2 }
		{ \Phi ( \mathrm{dist}_{ \mathbb{R} / \mathcal{L} \mathbb{Z} } ( s_1 , s_2 ) ) }
		\, d s_1 d s_2
		< 
		\infty
		\right. \right\}. 
	\]
	If $ \mathcal{E}_\Phi ( \vecf ) < \infty $, then 
	\[
		\vecf \in W_\Phi \cap W^{ 1, \infty } ( \mathbb{R} / \mathcal{L} \mathbb{Z} ) ,
		\quad
		\vecf \mbox{ is bi-Lipschitz}. 
	\]
	\item[{\rm (A.4)}]
	We set
	\[
		\Lambda (x) = - \frac 1x \int_x^\infty \frac { dt } { \Phi (t) }. 
	\]
	If $ \vecf  \in W_\Phi \cap W^{ 1, \infty } ( \mathbb{R} / \mathcal{L} \mathbb{Z} )  $ is bi-Lipschitz and 
	$ \| \vecf^\prime \|_{ \mathbb{R}^n } \equiv 1 $ almost everywhere,
	then it holds that
	\begin{align*}
		&
		\lim_{ \epsilon \to + 0 }
		\epsilon
		\int_{ \mathbb{R} / \mathcal{L} \mathbb{Z} }
		\left( \Lambda ( \| \vecf ( s_1 ) - \vecf ( s_1 + \epsilon ) \|_{ \mathbb{R}^n } )
		- \Lambda ( \epsilon ) \right)
		d s_1
		= 0 ,
		\\
		&
		\lim_{ \epsilon \to + 0 }
		\int_{ \mathbb{R} / \mathcal{L} \mathbb{Z} }
		\Lambda ( \| \vecf ( s_1 ) - \vecf ( s_1 + \epsilon ) \|_{ \mathbb{R}^n } )
		\int_{ s_1 }^{ s_1 + \epsilon }
		\| \vecf^\prime ( s_1 ) - \vecf^\prime ( s_2 ) \|_{ \mathbb{R}^n }^2
		 d s_2 d s_1
		 = 0. 
	\end{align*}
\item[{\rm (A.5)}]
	\begin{itemize}
	\item[{\rm (a)}]
	For any $ \lambda \in ( 0,1 ) $ and any $ x \in \left( 0 , \frac { \mathcal{L} } 2 \right] $, there is a constant $ C ( \lambda , \mathcal{L} ) > 0 $ such that $ \Phi ( \lambda x ) \geqq C( \lambda , \mathcal{L} ) \Phi ( x ) $, 
	\item[{\rm (b)}]
	$ \displaystyle{ \inf_{ x \in \left( 0 , \frac { \mathcal{L} } 2 \right] }
\left( \frac 1 { \Phi (x) } + \Lambda (x) \right) \geqq 0 } $.
	\end{itemize}
\end{itemize}
\par
In the case of O'Hara $ \Phi (t) = t^\alpha $,
the above assumptions hold if $ \alpha \in [ 2 , 3 ) $ where the functional $ \mathcal{E}_{ t^\alpha } = \mathcal{E}_{( \alpha , 1 ) } $ performs well as a knot energy (\cite{OH2}).
In particular,
see \cite{Blatt2} for (A.3).
\par
Let $ \psi ( s_1 , s_2 ) $ be the angle between $ \vectau ( s_1 ) $ and $ \vectau ( s_2 ) $.
We set
\[
	\Theta_{ \Phi } (t) = 
	\frac 12 \left(
	1
	+
	\Phi ( t ) \Lambda ( t )
	\right)
\]
The following is an extension of the cosine formula for generalized O'Hara energies.
\begin{thm}
\begin{enumerate}
\item
Under (A.1)--(A.4),
a generalized cosine formula
\begin{align*}
	\mathcal{E}_\Phi ( \vecf )
	= & \
	\mathrm{p.v.} \iint
	\left\{
	\left(
	1
	-
	\Theta_\Phi ( \| \Delta \vecf \| )
	\right)
	( 1 - \cos \varphi ( s_1 , s_2 ) )
	\right.
	\\
	& \qquad \qquad \qquad
	\left.
	+ \,
	\Theta_\Phi ( \| \Delta \vecf \| )
	( 1 - \cos \psi ( s_1 , s_2 ) )
	\right\}
	\frac { d s_1 d s_2 } { \Phi ( \| \Delta \vecf \| ) }
	\\
	& \quad
	+ \,
	2 \mathcal{L} \int_{ \frac { \mathcal{L} } 2 }^\infty
	\frac { dx } { \Phi (x) }
\end{align*}
holds.
Here $ \mathrm{p.v.} \iint = \lim_{ \epsilon \to + 0 } \iint_{ | s_1 - s_2 | \geqq \epsilon } $ is the integration in the principal value sense.
\item
In addition to (A.1)--(A.4),
if we assume (A.5) (b),
then
\[
	\varphi_{{}_\Phi} ( s_1 , s_2 )
	=
	\arccos
	\left\{
	\left(
	1
	-
	\Theta_\Phi ( \| \Delta \vecf \| )
	\right)
	\cos \varphi ( s_1 , s_2 )
	+
	\Theta_\Phi ( \| \Delta \vecf \| )
	\cos \psi ( s_1 , s_2 )
	\right\}
\]
can be defined,
and it holds that
\[
	\mathcal{E}_\Phi ( \vecf )
	=
	\iint_{ ( \mathbb{R} / \mathcal{L} \mathbb{Z} )^2 }
	\frac { 1 - \cos \varphi_{{}_\Phi} ( s_1 , s_2 ) } { \Phi ( \| \Delta \vecf \| ) }
	\, d s_1 d s_2
	+
	2 \mathcal{L} \int_{ \frac { \mathcal{L} } 2 }^\infty
	\frac { dx } { \Phi (x) }
	.
\]
\end{enumerate}
\label{Theorem}
\end{thm}
\par
We do not need (A.5) (a).
Proof of this will be given in the next section.
\par
Abrams-Cantarella-Fu-Ghomi-Howard \cite{ACFGH} established a condition under which the right circle minimizes the knot energies.
Our theorem establishes another type of condition.
Assume (A.1)--(A.4) and (A.5) (b).
Set
\begin{align*}
	\mathcal{E}_{ \Phi , 3 } ( \vecf )
	= & \
	\iint_{ ( \mathbb{R} / \mathcal{L} \mathbb{Z} )^2 }
	\left(
	1
	-
	\Theta_\Phi ( \| \Delta \vecf \| )
	\right)
	( 1 - \cos \varphi ( s_1 , s_2 ) )
	\frac { d s_1 d s_2 } { \Phi ( \| \Delta \vecf \| ) }
	,
	\\
	\mathcal{E}_{ \Phi , 4 } ( \vecf )
	= & \
	\iint_{ ( \mathbb{R} / \mathcal{L} \mathbb{Z} )^2 }
	\Theta_\Phi ( \| \Delta \vecf \| )
	( 1 - \cos \psi ( s_1 , s_2 ) )
	\frac { d s_1 d s_2 } { \Phi ( \| \Delta \vecf \| ) }
	.
\end{align*}
Let $ C $ be a circle with a circumstance that is the same as the total length $ \mathcal{L} $ of $ \mathrm{Im} \vecf $.
As the conformal angle for $ C $ vanishes identically,
we have $ \mathcal{E}_{ \Phi , 3 } ( C ) = 0 $.
Since the energy density of $ \mathcal{E}_{ \Phi , 3 } $ is non-negative,
we have
\[
	\mathcal{E}_\Phi ( \vecf )
	=
	\mathcal{E}_{ \Phi , 3 } ( \vecf )
	+
	(
	\mathcal{E}_{ \Phi , 4 } ( \vecf )
	-
	\mathcal{E}_{ \Phi , 4 } ( C )
	)
	+
	\mathcal{E}_\Phi ( C )
	\geqq
	(
	\mathcal{E}_{ \Phi , 4 } ( \vecf )
	-
	\mathcal{E}_{ \Phi , 4 } ( C )
	)
	+
	\mathcal{E}_\Phi ( C )
	.
\]
Consequently we obtain
\begin{cor}
Assume (A.1)--(A.4),
and (A.5) (b),
If $ C $ is a minimizer of $ \mathcal{E}_{ \Phi , 4 } $ under the length-constraint,
then it minimizes $ \mathcal{E}_\Phi $ under the constraint.
\end{cor}
\par
In the case of O'Hara $ \Phi (t) = t^\alpha $,
the function $ \Theta_\Phi $ is a constant,
and the assertion of Theorem \ref{Theorem} becomes
\begin{align*}
	\mathcal{E}_{ t^\alpha } ( \vecf )
	= & \
	\mathcal{E}_{(\alpha,1)} ( \vecf )
	\\
	= & \
	\left\{ 1 - \frac { \alpha - 2 } { 2 ( \alpha - 1 ) } \right\}
	\iint_{ ( \mathbb{R} / \mathcal{L} \mathbb{Z} )^2 }
	\frac { 1 - \cos \varphi ( s_1 , s_2 ) } { \| \Delta \vecf \|^\alpha }
	\, d s_1 d s_2
	\\
	& \quad
	+ \,
	\frac { \alpha - 2 } { 2 ( \alpha - 1 ) }
	\iint_{ ( \mathbb{R} / \mathcal{L} \mathbb{Z} )^2 }
	\frac { 1 - \cos \psi ( s_1 , s_2 ) } { \| \Delta \vecf \|^\alpha }
	\, d s_1 d s_2
	+
	\frac { 2^\alpha } { ( \alpha - 1 ) \mathcal{L}^{ \alpha - 2 } }
	.
\end{align*}
This coincides with \pref{cosine_formula_of_Moebius_energy} when $ \alpha = 2 $.
Now consider the normalized O'Hara energy:
\begin{align*}
	&
	\mathcal{L}^{ \alpha - 2 } \mathcal{E}_{( \alpha ,1 )} ( \vecf )
	\\
	& \quad
	=
	\left\{ 1 - \frac { \alpha - 2 } { 2 ( \alpha - 1 ) } \right\}
	\iint_{ ( \mathbb{R} / \mathcal{L} \mathbb{Z} )^2 }
	\left( \frac { \mathcal{L} } { \| \Delta \vecf \| } \right)^{ \alpha - 2 }
	( 1 - \cos \varphi ( s_1 , s_2 ) )
	\frac { d s_1 d s_2 } { \| \Delta \vecf \|^2 }
	\\
	& \quad \qquad
	+ \,
	\frac { \alpha - 2 } { 2 ( \alpha - 1 ) }
	\iint_{ ( \mathbb{R} / \mathcal{L} \mathbb{Z} )^2 }
	\left( \frac { \mathcal{L} } { \| \Delta \vecf \| } \right)^{ \alpha - 2 }
	\frac { 1 - \cos \psi ( s_1 , s_2 ) } { \| \Delta \vecf \|^2 }
	\, d s_1 d s_2
	+
	\frac { 2^\alpha } { \alpha - 1 }
	.
\end{align*}
It is scale-invariant.
The quantities
\[
	1 - \cos \varphi ( s_1 , s_2 )
	,
	\quad
	\frac { d s_1 d s_2 } { \| \Delta \vecf \|^2 }
	,
	\quad
	\iint_{ ( \mathbb{R} / \mathcal{L} \mathbb{Z} )^2 }
	\frac { 1 - \cos \psi ( s_1 , s_2 ) } { \| \Delta \vecf \|^2 }
	\, d s_1 d s_2
\]
are M\"{o}bius invariant (for the last one,
see \cite{IshizekiNagasawaI,IshizekiNagasawaIII}),
but not
\[
	\left( \frac { \mathcal{L} } { \| \Delta \vecf \| } \right)^{ \alpha - 2 }
\]
Consequently,
we can see from our formula how far the energy is from the M\"{o}bius invariant property.
\section{Proof}
\par
Employing the argument in \cite{IshizekiNagasawaV},
we can prove Theorem \ref{Theorem} very simply.
In \cite{IshizekiNagasawaV},
we proved
\begin{align}
\label{p.v.version}
	\mathcal{E}_\Phi ( \vecf )
	= & \
	\mathrm{p.v.} \iint
	\left[
	\frac { \| \Delta \vectau \|^2 } { 2 \Phi ( \| \Delta \vecf \| ) }
	\right.
	\\
	& \qquad \qquad
	\left.
	+ \,
	\left(
	\frac 1 { \Phi ( \| \Delta \vecf \| ) }
	-
	\Lambda ( \| \Delta \vecf \| )
	\right)
	\right.
	\nonumber
	\\
	& \quad \qquad \qquad
	\left.
	\times
	\left\{
	\vectau ( s_1 ) \cdot \vectau ( s_2 )
	-
	\left( \vectau ( s_1 ) \cdot \frac { \Delta \vecf } { \| \Delta \vecf \| } \right)
	\left( \vectau ( s_2 ) \cdot \frac { \Delta \vecf } { \| \Delta \vecf \| } \right)
	\right\}
	\right]
	d s_1 d s_2
	\nonumber
	\\
	& \quad
	+ \,
	2 \mathcal{L} \int_{ \frac { \mathcal{L} } 2 }^\infty
	\frac { dx } { \Phi (x) }
	\nonumber
\end{align}
under (A.1)--(A.4).
We can see that the conformal angle is given as
\[
	\cos \varphi ( s_1 , s_2 )
	=
	-
	\vectau ( s_1 ) \cdot \vectau ( s_2 )
	+
	2
	\left( \vectau ( s_1 ) \cdot \frac { \Delta \vecf } { \| \Delta \vecf \| } \right)
	\left( \vectau ( s_2 ) \cdot \frac { \Delta \vecf } { \| \Delta \vecf \| } \right)
	.
\]
Hence we have
\begin{align*}
	&
	\vectau ( s_1 ) \cdot \vectau ( s_2 )
	-
	\left( \vectau ( s_1 ) \cdot \frac { \Delta \vecf } { \| \Delta \vecf \| } \right)
	\left( \vectau ( s_2 ) \cdot \frac { \Delta \vecf } { \| \Delta \vecf \| } \right)
	\\
	& \quad
	=
	- \frac 12
	\left(
	1
	-
	\vectau ( s_1 ) \cdot \vectau ( s_2 )
	\right)
	\\
	& \quad \qquad
	+ \,
	\frac 12
	\left\{
	1
	+
	\vectau ( s_1 ) \cdot \vectau ( s_2 )
	-
	2
	\left( \vectau ( s_1 ) \cdot \frac { \Delta \vecf } { \| \Delta \vecf \| } \right)
	\left( \vectau ( s_2 ) \cdot \frac { \Delta \vecf } { \| \Delta \vecf \| } \right)
	\right\}
	\\
	& \quad
	=
	- \frac 14 \| \Delta \vectau \|^2
	+
	\frac 12 ( 1 - \cos \varphi ( s_1 , s_2 ) )
	.
\end{align*}
Using
\[
	\| \Delta \vectau \|^2
	= 2 \left( 1 - \vectau ( s_1 ) \cdot \vectau ( s_2 ) \right)
	= 2 ( 1 - \cos \psi ( s_1 , s_2 ) )
	,
\]
we obtain the first assertion of Theorem \ref{Theorem}.
Indeed,
the integrand of \pref{p.v.version} is
\begin{align*}
	&
	\frac 12
	\left(
	\frac 1 { \Phi ( \| \Delta \vecf \| ) }
	-
	\Lambda ( \| \Delta \vecf \| )
	\right)
	( 1 - \cos \varphi ( s_1 , s_2 ) )
	\\
	& \qquad \qquad
	+ \,
	\frac 14
	\left(
	\frac 1 { \Phi ( \| \Delta \vecf \| ) }
	+
	\Lambda ( \| \Delta \vecf \| )
	\right)
	\| \Delta \vectau \|^2
	\\
	& \quad
	=
	\frac 1 { 2 \Phi ( \| \Delta \vecf \| ) }
	\left\{
	\left(
	1
	-
	\Phi ( \| \Delta \vecf \| ) \Lambda ( \| \Delta \vecf \| )
	\right)
	( 1 - \cos \varphi ( s_1 , s_2 ) )
	\right.
	\\
	& \qquad \qquad \qquad \qquad
	\left.
	+ \,
	\left(
	1
	+
	\Phi ( \| \Delta \vecf \| ) \Lambda ( \| \Delta \vecf \| )
	\right)
	( 1 - \cos \psi ( s_1 , s_2 ) )
	\right\}
	\\
	& \quad
	=
	\frac 1 { \Phi ( \| \Delta \vecf \| ) }
	\left\{
	\left(
	1
	-
	\Theta_\Phi ( \| \Delta \vecf \| )
	\right)
	( 1 - \cos \varphi ( s_1 , s_2 ) )
	\right.
	\\
	& \qquad \qquad \qquad \qquad
	\left.
	+ \,
	\Theta_\Phi ( \| \Delta \vecf \| )
	( 1 - \cos \psi ( s_1 , s_2 ) )
	\right\}
	.
\end{align*}
\par
If we also assume (A.5) (b),
then
\[
	 0 \leqq \Theta (x) \leqq \frac 12 .
\]
Therefore,
the integrand
\[
	\left(
	1
	-
	\Theta_\Phi ( \| \Delta \vecf \| )
	\right)
	( 1 - \cos \varphi ( s_1 , s_2 ) )
	+
	\Theta_\Phi ( \| \Delta \vecf \| )
	( 1 - \cos \psi ( s_1 , s_2 ) )
\]
is non-negative,
and
\[
	\left|
	\left(
	1
	-
	\Theta_\Phi ( \| \Delta \vecf \| )
	\right)
	\cos \varphi ( s_1 , s_2 ) )
	+
	\Theta_\Phi ( \| \Delta \vecf \| )
	\cos \psi ( s_1 , s_2 ) )
	\right|
	\leqq
	1
	.
\]
Hence,
the integration in the principal value sense becomes one in the usual $ L^1 $ sense.
Furthermore,
the angle $ \varphi_{{}_\Phi} $ is defined.
This shows the second assertion of Theorem \ref{Theorem}.

\end{document}